\documentclass[12pt,a4paper]{amsart}

\usepackage{amssymb}
\usepackage{amsmath}
\usepackage{amscd}
\usepackage{amsbsy}
\usepackage[matrix,arrow]{xy}
\usepackage{comment}
\usepackage{enumerate}

\DeclareMathOperator{\ord}{ord}

\newcommand{\Q}{{\mathbb Q}}
\newcommand{\Z}{{\mathbb Z}}

\newcommand{\fp}{\mathfrak{p}}
\newcommand{\fq}{\mathfrak{q}}

\newcommand{\eps}{\varepsilon}

\newtheorem{thm}{Theorem}[section]

\newtheorem{prop}[thm]{Proposition}

\theoremstyle{definition}

\theoremstyle{remark}

\begin{document}

\title{The Generalised Fermat Equation $x^2+y^3=z^{15}$}
\author{Samir Siksek}
\address{Mathematics Institute\\
	University of Warwick\\
	Coventry\\
	CV4 7AL \\
	United Kingdom}

\email{s.siksek@warwick.ac.uk}

\author{Michael Stoll}
\address{Mathematisches Institut,
         Universit\"at Bayreuth,
         95440 Bayreuth, Germany.}
\email{Michael.Stoll@uni-bayreuth.de}

\date{12 January 2014}
\thanks{The first-named author is supported by an EPSRC Leadership Fellowship.}


\begin {abstract}
  We determine the set of primitive integral solutions to the
  generalised Fermat equation $x^2 + y^3 = z^{15}$. As expected,
  the only solutions are the trivial ones with $xyz = 0$ and
  the non-trivial pair $(x,y,z) = (\pm 3, -2, 1)$.
\end {abstract}

\maketitle


\section{Introduction}

Let $p$, $q$, $r \in \Z_{\geq 2}$. The equation
\begin{equation}\label{eqn:FCgen}
x^p+y^q=z^r
\end{equation}
is known as the \emph{generalised Fermat equation}
(or the \emph{Fermat-Catalan} equation) with signature $(p,q,r)$.
As in Fermat's Last Theorem, one is interested in non-trivial
primitive integer solutions.
An integer solution $(x,y,z)$ is said to be {\em non-trivial} if
$xyz \neq 0$, and {\em primitive} if $x$, $y$, $z$ are coprime.
Let $\chi=p^{-1}+q^{-1}+r^{-1}$. The parametrisation
of non-trivial primitive integer solutions for $(p,q,r)$ with
$\chi \geq 1$ has been completed~\cite{Ed}.
The generalised Fermat Conjecture~\cite{Da97,DG}
is concerned with the case $\chi<1$.
It states that the only non-trivial primitive solutions to
\eqref{eqn:FCgen} with $\chi<1$ are (up to sign and permutation)
\begin{gather*}
1+2^3 = 3^2, \quad 2^5+7^2 = 3^4, \quad 7^3+13^2 = 2^9, \quad
2^7+17^3 = 71^2, \\
3^5+11^4 = 122^2, \quad 17^7+76271^3 = 21063928^2, \quad
1414^3+2213459^2 = 65^7, \\
9262^3+15312283^2 = 113^7, \quad
43^8+96222^3 = 30042907^2, \quad 33^8+1549034^2 = 15613^3.
\end{gather*}
The generalised Fermat Conjecture
has been established for many signatures $(p,q,r)$,
including for several infinite families of signatures,
see for example~\cite{SikSto} for a short overview,
and \cite[Chapter 17]{Cohen} for a relatively recent survey.

There is an abundance of solutions for the generalised Fermat equation
with signatures $(2,3,n)$, and so this subfamily is particularly interesting.
The condition $\chi<1$ within this family coincides with $n \ge 7$. 
The cases $n=7$, $8$, $9$ are solved respectively in
\cite{PSS}, \cite{Bruin}, \cite{Br3}. The case $n=10$ is solved
independently in \cite{Brown} and~\cite{chabnf}.
Every $n>5$ is divisible by 
$6$, $8$, $9$, $10$, $15$, $25$ or an
odd prime $p \ge 7$. Given the results for $6 \le n \le 10$, it would be
sufficient for a complete resolution of this subfamily to deal with
exponents $n=15$, $25$ and prime exponents $n \ge 11$.
In this note, we deal with the case $n = 15$. Our result is as follows.

\begin{thm}\label{thm:main}
  The only primitive integer solutions to the equation
  \begin{equation}\label{eqn:main}
    x^2 + y^3 = z^{15}
  \end{equation}
  are the trivial solutions $(\pm 1, -1, 0)$, $(\pm 1, 0, 1)$, $(0, 1, 1)$,
  $(0, -1, -1)$ and the non-trivial solutions $(\pm 3, -2, 1)$.
\end{thm}

There are two fairly obvious avenues for attacking equation~\eqref{eqn:main}.
One is to use Edwards' parametrisations \cite{Ed} of the solutions to $x^2 + y^3 = z^5$.
In these parametrisations, $z$ is given by a binary form of degree~$12$,
so the requirement that $z$ be a cube results in a number of superelliptic
curves of genus~$10$ of the form $w^3 = f(u,v)$ with $\deg f = 12$.
To these curves (and if necessary, to their Jacobians), one can apply descent
techniques analogous to those discussed in~\cite{Mourao},~\cite{PS},~\cite{SikSto}.

However, we decided to follow the other route, which is based on the parametrisations
of the solutions to $x^2 + y^3 = z^3$. In this case, the equations obtained
are of the form $w^5 = f(u,v)$ with a binary quartic form~$f$. A further
descent step reduces the problem to that of finding rational points on a
number of hyperelliptic curves of genus~$2$ or~$4$.


\section{Parametrisation of Solutions to $x^2=y^3+z^3$}

We recall the following result \cite[Theorem 14.3.1]{Cohen}.

\begin{thm} \label{thm:233}
  Let $(x,y,z)$ be a triple of coprime integers such that $x^2 = y^3 + z^3$.
  Then, up to possibly interchanging $y$ and~$z$, there are coprime integers
  $s$ and~$t$ such that one of the following sets of relations hold.
  \begin{enumerate}[\upshape (P1)]
    \renewcommand{\arraystretch}{1.2}
    \addtolength{\itemsep}{3pt}
    \item $\displaystyle \left\{ \begin{array}{r@{{}={}}l}
            x & \pm (s^2 - 2 s t - 2 t^2)(s^4 + 2 s^3 t + 6 s^2 t^2 - 4 s t^3 + 4 t^4) \\
            y & s (s + 2t) (s^2 - 2 s t + 4 t^2) \\
            z & -4 t (s - t) (s^2 + s t + t^2)
                     \end{array} \right.
          $ \\
          with $s \not\equiv 0 \bmod 2$ and $s \not\equiv t \bmod 3$.
    \item $\displaystyle \left\{ \begin{array}{r@{{}={}}l}
            x & 3 (s - t)(s + t)(s^4 + 2 s^3 t + 6 s^2 t^2 + 2 s t^3 + t^4) \\
            y & s^4 - 4 s^3 t - 6 s^2 t^2 - 4 s t^3 + t^4 \\
            z & 2 (s^4 + 2 s^3 t + 2 s t^3 + t^4)
                     \end{array} \right.
          $ \\
          with $s \not\equiv t \bmod 2$ and $s \not\equiv t \bmod 3$.
    \item $\displaystyle \left\{ \begin{array}{r@{{}={}}l}
            x &= 6 s t (3 s^4 + t^4) \\
            y &= -3 s^4 + 6 s^2 t^2 + t^4 \\
            z &= 3 s^4 + 6 s^2 t^2 - t^4
                     \end{array} \right.
          $ \\
          with $s \not\equiv t \bmod 2$ and $t \not\equiv 0 \bmod 3$.
  \end{enumerate}
\end{thm}
As far as we know, this parametrisation was first obtained by Mordell \cite[page 235]{Mordell},
but see the discussion in \cite[Appendix B.1]{Ed}.

As a corollary, we see that each primitive integral solution to $x^2 + y^3 = z^{15}$
gives rise to a solution of an equation
\begin{equation}\label{eqn:ufj}
 u^5 = f_i(s, t) 
\end{equation}
for some $i=1,2,\ldots,6$,
where $f_1, \ldots, f_6$ are the six binary quartic forms giving the values of $y$
and~$z$ in Theorem~\ref{thm:233}, with $s$, $t$ coprime and satisfying the relevant conditions
modulo $2$ and $3$.
For future reference, we define the $f_i$ explicitly as follows.
\begin{align*}
  f_1(s, t) &= s (s + 2t) (s^2 - 2 s t + 4 t^2) \\
  f_2(s, t) &= -4 t (s - t) (s^2 + s t + t^2) \\
  f_3(s, t) &= s^4 - 4 s^3 t - 6 s^2 t^2 - 4 s t^3 + t^4 \\
  f_4(s, t) &= 2 (s^4 + 2 s^3 t + 2 s t^3 + t^4) \\
  f_5(s, t) &= -3 s^4 + 6 s^2 t^2 + t^4 \\
  f_6(s, t) &= 3 s^4 + 6 s^2 t^2 - t^4
\end{align*}
We note that we can eliminate $f_4$ immediately: here $s$ and~$t$ must
be of different parity, which implies that $s^4 + 2 s^3 t + 2 s t^3 + t^4$
is odd, so $f_4(s, t)$ is divisible by~$2$, but not by~$4$ and hence cannot be
a fifth power. The known  solutions of equation~\eqref{eqn:main}
give rise to solutions of $u^5 = f_i(s, t)$ for $i \in \{1,2,3,5,6\}$, which
means that it is not possible to rule out any of the other quartic forms
on the basis of local considerations. In what follows we
shall deal separately with equation \eqref{eqn:ufj} with $i=1,2,3,5,6$
and $s$, $t$ satisfying the relevant conditions. We shall refer
to these as Cases $1$, $2$, $3$, $5$ and $6$, which are respectively
resolved in Sections~\ref{sec:case1}, \ref{sec:case2}, \ref{sec:case3},
\ref{sec:case5}, \ref{sec:case6}.
All computations mentioned below were carried out using the computer
algebra package {\tt MAGMA} \cite{MAGMA}.

\section{Case 1}\label{sec:case1}

We want to solve the equation
\[ u^5 = s (s + 2t) (s^2 - 2 s t + 4 t^2) \]
in integers $s,t,u$ with $s$ and~$t$ coprime, $s \not\equiv 0 \bmod 2$
and $s \not\equiv t \bmod 3$. These congruence and coprimality conditions imply
that the three factors on the right hand side are coprime in pairs.
This in turn implies that there are integers $w_1, w_2, w_3$, coprime in pairs,
such that
\begin{equation}\label{eqn:case1}
 s = w_1^5, \quad s + 2t = w_2^5, \quad s^2 - 2 s t + 4 t^2 = w_3^5 \,.
\end{equation} 
If we substitute the first of these equations into the last, we obtain
\[ \bigl(2 t - \tfrac{1}{2} w_1^5\bigr)^2 + \tfrac{3}{4} w_1^{10}
     = w_1^{10} - 2 t w_1^5 + 4 t^2 = w_3^5 \,,
\]
which, upon setting $X = w_3/w_1^2$ and $Y = 4t/w_1^5 - 1$, becomes
\[ Y^2 = 4 X^5 - 3 \,. \]
This is a curve of genus~$2$, to which by now standard methods can be applied.
A $2$-descent on its Jacobian as described in~\cite{Stoll} gives an upper
bound of~$1$ for the Mordell-Weil rank. An argument using the explicit theory of
heights as developed in~\cite{Stollh1,Stollh2} shows that the point on the Jacobian
given by the class of the divisor
\[ (i, 2i+1) + (-i, -2i+1) - 2 \infty \]
(where $i^2 = -1$) generates the Mordell-Weil group (it is easy to check that
the torsion subgroup is trivial). Finally, a computation using Chabauty's
method and the Mordell-Weil sieve as explained in~\cite{BS3} shows that
the rational points on the curve are
\[ \infty, \quad (1, 1), \quad (1, -1) \,. \]
The first of these leads to $s = 0$, which is excluded ($s$ must be odd).
The second leads to the contradiction $2 w_1^5 = w_2^5$, whereas the last
gives the trivial solutions $(x,y,z) = (\pm 1, 0, 1)$. We have shown
the following.

\begin{prop} \label{prop:f1}
  The only solutions in integers $(s,t,u)$ with $s$ and~$t$ coprime and satisfying
  $s \not\equiv 0 \bmod 2$ and $s \not\equiv t \bmod 3$ to the equation
  \[ u^5 = s (s + 2t) (s^2 - 2 s t + 4 t^2) \]
  are $(s,t,u) = (\pm 1, 0, 1)$. 
This yields the trivial
solutions $(x,y,z)=(\pm 1, 0, 1)$ to \eqref{eqn:main}.
\end{prop}

\section{Case 2}\label{sec:case2}

Now we consider $f_2$, so we want to solve
\[ u^5 = -4 t (s - t) (s^2 + s t + t^2) \]
in integers with $s$ and~$t$ coprime, $s \not\equiv 0 \bmod 2$,
$s \not\equiv t \bmod 3$.
As before, these conditions imply that the three non-constant factors on the
right-hand side are coprime in pairs, and the last factor is odd.
Therefore there are integers $w_1, w_2, w_3$, coprime in pairs, such that
\begin{equation}\label{eqn:case2.1}
 t= w_1^5, \quad s-t=8 w_2^5, \quad s^2 + s t + t^2 = w_3^5 \, , 
\end{equation}
or
\begin{equation}\label{eqn:case2.2}
 t = 8 w_1^5, \quad s - t = w_2^5, \quad s^2 + s t + t^2 = w_3^5 \,. 
\end{equation}
We deal with the possibilities \eqref{eqn:case2.1} and \eqref{eqn:case2.2}
as we have dealt with \eqref{eqn:case1}. If \eqref{eqn:case2.1} holds then
\[
(2s + t)^2 + 3 t^2  = 4(s^2 + s t+ t^2) = 4 w_3^5 
\]
which shows that $(w_3/w_1^2, (2s+w_1^5)/ w_1^5 )$ is a point on the curve $Y^2=4 X^5-3$
whose points we have already determined; we obtain $(s,t,u)= (\pm 1, 0 , 0)$,
and $(x,y,z)=(\pm 1, 1, 0)$.

If \eqref{eqn:case2.2} holds then
\[ (s + 4 w_1^5)^2 + 48 w_1^{10}
     = \bigl(s + \tfrac{1}{2}t\bigr)^2 + \tfrac{3}{4} t^2 = w_3^5 \,,
\]
leading to the genus~$2$ curve
\[ Y^2 = X^5 - 48 \]
(with $X = w_3/w_1^2$ and $Y = s/w_1^5 + 4$). As before, the Mordell-Weil
group is infinite cyclic, this time generated by the class of
\[ (6 + 2 \sqrt{2}, 124 + 76 \sqrt{2}) + (6 - 2 \sqrt{2}, 124 - 76 \sqrt{2}) - 2 \infty \, ,\]
and a `Chabauty plus Mordell-Weil sieve' computation shows that this curve
has the point at infinity as its only rational point. The solutions corresponding
to this have $w_1 = 0$, so $t = 0$ and $s = \pm 1$, $u = 0$. This gives the
trivial solutions $(x,y,z) = (\pm 1, 1, 0)$ of the original equation~\eqref{eqn:main}.
We have shown:

\begin{prop} \label{prop:f2}
  The only solutions in integers $(s,t,u)$ with $s$ and~$t$ coprime and satisfying
  $s \not\equiv 0 \bmod 2$ and $s \not\equiv t \bmod 3$ to the equation
  \[ u^5 = -4 t (s - t) (s^2 + s t + t^2) \]
  are $(s,t,u) = (\pm 1, 0, 0)$.
This yields the trivial
solutions $(x,y,z)=(\pm 1, 1, 0)$ to \eqref{eqn:main}.
\end{prop}

\section{Case 3}\label{sec:case3}

The remaining three forms $f_3$, $f_5$ and~$f_6$ are all irreducible.
However, they all split over $\Q(\sqrt{3})$ into two quadratic factors
that can be written as linear combinations of two squares of linear forms
over~$\Q$. The two factors must be again coprime (in $\Z[\sqrt{3}]$) and
are therefore fifth powers up to a power of the fundamental unit $\eps = \sqrt{3} - 2$.
We then express two squares of linear forms as binary quintic forms;
taking their product results in a hyperelliptic curve of genus~$4$.
In this section we carry this out for~$f_3$. We have
\[
\begin{split} 
f_3(s,t) & = s^4 - 4 s^3 t - 6 s^2 t^2 - 4 s t^3 + t^4\\
         & = \bigl(s^2 - 2(1+\sqrt{3}) s t + t^2\bigr)\bigl(s^2 - 2(1-\sqrt{3}) s t + t^2\bigr) \, .
\end{split}
\]
The conditions are that $s$ and~$t$ are coprime with $s \not\equiv t \bmod 2$
and $s \not\equiv t \bmod 3$, which imply that each factor is coprime to~$6$.
Since their resultant is~$48$, it follows that they are coprime. So there
is $j \in \{-2,-1,0,1,2\}$ and there are integers $v, w$ such that
\[
s^2 - 2(1-\sqrt{3}) s t + t^2 
     = \eps^j (v + w \sqrt{3})^5
     = g_j(v, w) + h_j(v, w) \sqrt{3}
\]
with binary quintic forms $g_j$, $h_j$ with rational integral coefficients.
Thus
\begin{equation}\label{eqn:gh}
s^2 - 2 s t + t^2 = g_j(v,w), \qquad 2 s t = h_j(v,w) \, .
\end{equation}
Explicitly, we have
\begin{align*}
  g_{-2}(v,w) &= 7 v^5 + 60 v^4 w + 210 v^3 w^2 + 360 v^2 w^3 + 315 v w^4 + 108 w^5 \\
  g_{-1}(v,w) &= -(2 v^5 + 15 v^4 w + 60 v^3 w^2 + 90 v^2 w^3 + 90 v w^4 + 27 w^5) \\
  g_0(v,w)    &= v (v^4 + 30 v^2 w^2 + 45 w^4) \\
  g_1(v,w)    &= -2 v^5 + 15 v^4 w - 60 v^3 w^2 + 90 v^2 w^3 - 90 v w^4 + 27 w^5 = g_{-1}(v,-w) \\
  g_2(v,w)    &= 7 v^5 - 60 v^4 w + 210 v^3 w^2 - 360 v^2 w^3 + 315 v w^4 - 108 w^5 = g_{-2}(v,-w)
\end{align*}
and
\begin{align*}
  h_{-2}(v,w) &= 4 v^5 + 35 v^4 w + 120 v^3 w^2 + 210 v^2 w^3 + 180 v w^4 + 63 w^5 \\
  h_{-1}(v,w) &= -v^5 - 10 v^4 w - 30 v^3 w^2 - 60 v^2 w^3 - 45 v w^4 - 18 w^5 \\
  h_0(v,w)    &= (5 v^4 + 30 v^2 w^2 + 9 w^4) w \\
  h_1(v,w)    &= v^5 - 10 v^4 w + 30 v^3 w^2 - 60 v^2 w^3 + 45 v w^4 - 18 w^5 = h_{-1}(-v,w) \\
  h_2(v,w)    &= -4 v^5 + 35 v^4 w - 120 v^3 w^2 + 210 v^2 w^3 - 180 v w^4 + 63 w^5 = h_{-2}(-v,w) \, .
\end{align*}
From \eqref{eqn:gh},
\[ (s-t)^2 = g_j(v, w) \quad\text{and}\quad
   (s+t)^2 = g_j(v, w) + 2 h_j(v, w)
\]
so, setting $Y = (s^2 - t^2)/w^5$, $X = v/w$, this gives the hyperelliptic curve
\[ C_{3,j} \colon Y^2 = g_j(X,1)\bigl(g_j(X,1) + 2 h_j(X,1)\bigr) \]
of genus~$4$.
The irreducible quintic factors occurring on the right hand side of these equations
all have a root in the same number field $L = \Q(\theta)$ where
\begin{equation}\label{eqn:theta}
 \theta^5 - 5 \theta^3 + 5 \theta - 4 = 0\,. 
\end{equation}
A $2$-cover descent
as described in~\cite{BS2} on the first two and the last curves ($j=-2$, $-1$, $2$)
proves that they do not
possess any rational points. On the other hand, the two remaining curves do have
some rational points, so we need to put in some more work. We consider $C_{3,1}$ first.
A partial $2$-descent as in~\cite{SikSto} using the factorisation over~$\Q$ shows that
if $(X,Y)$ is a rational point on~$C_{3,1}$, then there is a rational point $(X,\tilde{Y})$
(with the same $X$-coordinate!) on either
\[ {\tilde{Y}}^2 = 2 X^5 - 15 X^4 + 60 X^3 - 90 X^2 + 90 X - 27 \]
or
\[ 5 {\tilde{Y}}^2 = 2 X^5 - 15 X^4 + 60 X^3 - 90 X^2 + 90 X - 27 \,. \]
These are curves of genus~$2$ again. For the first curve, we find in a similar way
as we did when considering $f_1$ and~$f_2$ that its Jacobian has Mordell-Weil rank~$1$,
and a `Chabauty plus MWS' computation shows that the rational points are the point
at infinity and two points with $X = 3$; the two points with $X=3$
do not give rise to points on~$C_{3,1}$.
For the second curve, the rank is again~$1$, and its rational points are the point
at infinity and two points with $X = 1$; again the points with
$X=1$ do not lead to points on~$C_{3,1}$.
Thus
\[
C_{3,1}(\Q)=\{\infty\} \, .
\]
As $X=v/w$, we see that $w=0$ and so $v=\pm 1$. This does not lead to any solutions
of~\eqref{eqn:main}.

In principle, one could try the same approach with~$C_{3,0}$,
which has equation
\begin{equation}\label{eqn:C30}
Y^2=X (X^4+30 X^2+ 45)(X^5 + 10 X^4 + 30 X^3 + 60 X^2 + 45 X + 18)
\end{equation}
Here one obtains the
two genus~$2$ curves
\[ {\tilde{Y}}^2 = X (X^5 + 10 X^4 + 30 X^3 + 60 X^2 + 45 X + 18) \]
and
\[ 5 {\tilde{Y}}^2 = X (X^5 + 10 X^4 + 30 X^3 + 60 X^2 + 45 X + 18) \,. \]
The second of these has Jacobian with trivial Mordell-Weil group and therefore only
the rational point $(0,0)$. However, the first curve has Jacobian of Mordell-Weil rank~$2$; therefore Chabauty is not applicable.
Note that the genus~$1$ curve $Y^2 = X^4 + 30 X^2 + 45$ (which one could also consider)
is an elliptic curve of positive rank and so is of little use for our purposes.
Instead it is more convenient to work
directly with the Jacobian $J_{3,0}$ of~$C_{3,0}$.
A $2$-descent on $J_{3,0}$ shows that its $2$-Selmer rank is~$2$,
and a search for reasonably small rational points on the 2-coverings of~$J_{3,0}$
corresponding to Selmer group elements
leads to the following points on the Jacobian.
\begin{gather*}
D_1=(X^4 + 30X^2 + 45,\; 0) \, , \qquad
 D_2=(0,0)- \infty_{+}\, , \\
D_3=  \biggl(
X^4 + \frac{8}{5} X^3 + 14X^2 + \frac{72}{5} X + \frac{81}{5} \, , \; \,
    -\frac{264}{25} X^3 + 32 X^2 + \frac{1304}{25} X + \frac{1152}{25}
\biggr)\, .
\end{gather*}
The points $D_1$ and $D_3$ are given in their Mumford representation:
$D=(f(X),g(X))$
means that $D=D^\prime-(d/2)(\infty_+
+\infty_{-})$ where $D^\prime$ is the effective degree $d$ divisor cut out
on the affine model \eqref{eqn:C30} by the 
simultaneous pair of equations $f(X)=Y-g(X)=0$.
The points $D_1$, $D_2$, $D_3$ generate a subgroup of $J(\Q)$
isomorphic to $\Z/2\Z \times \Z^2$. 
Hence the Mordell--Weil
rank is precisely $2$, and as the genus of $C_{3,0}$ is $4$,
Chabauty's method is applicable. A \lq Chabauty plus MWS\rq\
computation shows that the only rational points on $C_{3,0}$
are the $(0,0)$ and the two points at infinity.
Hence $X=0$ or $\infty$ and so $v=0$ or $w=0$.
By \eqref{eqn:gh}, $v=0$ gives $s=t$ which contradicts
$s \not\equiv t \bmod 2$, and so does not lead to a solution of \eqref{eqn:main}.
However $w = 0$ implies (again from \eqref{eqn:gh}) that
$s = 0$ or $t = 0$, corresponding to the non-trivial solutions $(x,y,z) = (\pm 3, -2, 1)$
of equation~\eqref{eqn:main}. We have shown:

\begin{prop} \label{prop:f3}
  The only solutions in integers $(s,t,u)$ with $s$ and~$t$ coprime and satisfying
  $s \not\equiv t \bmod 2$ and $s \not\equiv t \bmod 3$ to the equation
  \[ u^5 = s^4 - 4 s^3 t - 6 s^2 t^2 - 4 s t^3 + t^4 \]
  are $(s,t,u) = (\pm 1, 0, 1)$ and $(0, \pm 1, 1)$.
These yield the non-trivial
solutions $(x,y,z)=(\pm 3, -2, 1)$ to \eqref{eqn:main}.
\end{prop}

\section{Case 5}\label{sec:case5}
In this section we apply the method of the previous section to $f_5$.
We have
\[ f_5(s,t) = -3 s^4 + 6 s^2 t^2 + t^4
            = \bigl((3 + 2 \sqrt{3}) s^2 + t^2\bigr)\bigl((3 - 2 \sqrt{3}) s^2 + t^2\bigr) \,.
\]
The conditions on $s$ and~$t$ are coprimality and $s \not\equiv t \bmod 2$,
$t \not\equiv 0 \bmod 3$. They again imply that the two factors are coprime.
Writing
\[ (3 + 2 \sqrt{3}) s^2 + t^2 = \eps^j (v + w \sqrt{3})^5
                              = g_j(v, w) + h_j(v, w) \sqrt{3} \,,
\]
we find that
\begin{equation}\label{eqn:50}
 2 s^2 = h_j(v, w) \quad\text{and}\quad 2 t^2 = 2 g_j(v,w) - 3 h_j(v,w) 
\end{equation}
so that
\[ (2 s t)^2 = h_j(v, w) \bigl(2 g_j(v,w) - 3 h_j(v,w)\bigr) \,. \]
This defines again five hyperelliptic curves~$C_{5,j}$ of genus~$4$. It can be checked
that (taking $j$ mod~$5$) $C_{3,j} \cong C_{5,j-1}$. 
We know from the previous section that $C_{3,j}(\Q)=\emptyset$
for $j=-2$, $-1$ and $2$. Thus $C_{5,j}(\Q)=\emptyset$ for 
$j=2$, $-2$ and $1$.

The models given above for $C_{5,0}$ and $C_{3,1}$ are in fact identical, 
so by the results of the previous section $C_{5,0}(\Q)=\{\infty\}$.
The point at infinity on $C_{5,0}$ gives $w = 0$,
which this time gives $s = 0$ and $2 t^2 = 2 v^5$, hence the solutions
$(s,t) = (0, \pm 1)$, corresponding to the solution $(x,y,z) = (0,1,1)$ to \eqref{eqn:main}. 

The models given above for $C_{5,-1}$ and $C_{3,0}$ are also identical,
so we know again by the results of the previous section that
$C_{5,-1} = \{\infty_+, \infty_-, (0,0)\}$.
The points at infinity correspond to $w = 0$, which leads to $s = \pm t$, a contradiction.
The other point gives $v = 0$, which leads to $t = 0$, another contradiction.

We obtain the following.

\begin{prop} \label{prop:f5}
  The only solutions in integers $(s,t,u)$ with $s$ and~$t$ coprime and satisfying
  $s \not\equiv t \bmod 2$ and $t \not\equiv 0 \bmod 3$ to the equation
  \[ u^5 = -3 s^4 + 6 s^2 t^2 + t^4 \]
  are $(s,t,u) = (0, \pm 1, 1)$.
This yields the trivial solution $(x,y,z)=(0,1,1)$ to
\eqref{eqn:main}.
\end{prop}

\section{Case 6}\label{sec:case6}

To complete the proof of Theorem~\ref{thm:main} it remains to deal with $f_6$,
which we can write as
\[ f_6(s,t) = 3 s^4 + 6 s^2 t^2 - t^4
            = - \bigl(t^2 - 3 s^2 +2 s^2 \sqrt{3}\bigr)\bigl(t^2-3 s^2 -2 s^2 \sqrt{3}\bigr) \,.
\]
The conditions on $s$ and~$t$ are coprimality, $s \not\equiv t \bmod 2$ and
$t \not\equiv 0 \bmod 3$. They again imply that the two factors are coprime.
Writing
\[ t^2 - 3 s^2 +2 s^2 \sqrt{3} = \eps^j (v + w \sqrt{3})^5
                              = g_j(v, w) + h_j(v, w) \sqrt{3} \,,
\]
we find that
\begin{equation}\label{eqn:case6}
 2 s^2 = h_j(v, w) \quad\text{and}\quad 2 t^2 = 2 g_j(v,w) + 3 h_j(v,w) 
\end{equation}
so that
\[ (2 s t)^2 = h_j(v, w) \bigl(2 g_j(v,w) + 3 h_j(v,w)\bigr) \,. \]
We therefore obtain five curves~$C_{6,j}$, which
are easily seen to be the quadratic twists by~$-1$ of the curves~$C_{5,-j}$.
We deal with $C_{6,-2}$ last as it requires slightly more
delicate arguments.
We can rule out~$C_{6,-1}$ via a 2-cover descent. Partial $2$-descent on~$C_{6,0}$
shows that rational points give rise to rational points with the same $X$-coordinate
on the genus~$2$ curve
\[ {\tilde{Y}}^2 = 10 X^5 + 75 X^4 + 300 X^3 + 450 X^2 + 450 X + 135 \,, \]
which has Jacobian with trivial Mordell-Weil group. This implies that the only
rational point on~$C_{6,0}$ is the point at infinity. This gives the solutions
$(s,t) = (0, \pm 1)$, corresponding to $(x,y,z) = (0,-1,-1)$.

We consider the rational points on $C_{6,1}$. A partial $2$-descent shows
that rational points give rise to rational points with the same $X$-coordinate
on the genus~$1$ curve
\[ 5{\tilde{Y}}^2 =  X^4 + 30 X^2 + 45. \]
This curve has Mordell-Weil rank $1$ and so infinitely many
rational points. However, considering the equation $3$-adically,
we see that $3$-adic points must satisfy $\ord_3(X) \ge 1$.
As $X=v/w$, we see that $3 \mid v$. From \eqref{eqn:case6} we
have
\[ 2 t^2 =-v (v^4+30 v^2 w^2+45 w^4) \]
which contradicts $t \not \equiv 0 \pmod{3}$. Thus the rational points
on $C_{6,1}$ do not give rise to solutions to \eqref{eqn:main}.

Let us deal with $j=2$. It turns out to be more convenient to
work directly with \eqref{eqn:case6} than with the curve $C_{6,2}$.
Explicitly, \eqref{eqn:case6} is the following pair of equations
\[
\begin{split}
2 s^2 &= -4v^5 + 35v^4 w - 120 v^3 w^2 + 210 v^2 w^3 - 180 v w^4 + 63 w^5\\
2 t^2 &= 2 v^5 - 15 v^4 w + 60 v^3 w^2 - 90 v^2 w^3 + 90 v w^4 - 27 w^5 \, .
\end{split}
\]
Since $s \not\equiv t \mod 2$, we see that $v$, $w$ cannot both be
even. In fact, by listing all the possibilities for $v$, $w$ modulo $8$,
we see that $v$ must be odd and $w$ must be even. Now we take a
linear combination that eliminates the $v^5$ term:
\[ 2 (s^2+2 t^2)= w (5 v^4+30 v^2 w^2 +9 w^4). \] 
Recall that $s$, $t$ are coprime. Any odd prime $p$ dividing the
left-hand side must satisfy $\left(-2/p\right)=1$ 
and so $p \equiv 1$ or $3 \bmod {8}$. Thus, 
$5 v^4+30 v^2 w^2+9 w^4 \equiv 1$ or $3 \bmod 8$. 
As $v$ is odd and $w$ is even we see that
$5 v^4+30 v^2 w^2 +9 w^4 \equiv 5 \bmod 8$,
giving a contradiction.

\bigskip

The curve $C_{6,-2}$ has equation
\begin{multline*}
Y^2 = (4 X^5 + 35 X^4 + 120 X^3 + 210 X^2 + 180 X + 63)  \\
         \times (26 X^5 + 225 X^4 + 780 X^3 + 1350 X^2 + 1170 X + 405).
\end{multline*}
The two irreducible quintic factors on the right-hand side 
each acquire a root over $L=\Q(\theta)$ where
$\theta$ is given by \eqref{eqn:theta}. These roots are 
respectively
\[
\phi_1=\frac{\theta^4 - 5\theta^2 - 4\theta - 3}{4},
\qquad
\phi_2=\frac{7\theta^4 - 2\theta^3 - 27\theta^2 + 4\theta - 33}{26} \, .
\]
Let
\[
\mu_1=2 \theta^4 + 2\theta^3 - 3\theta^2 + \theta + 1, \qquad
\mu_2=\frac{18 \theta^4 + 19 \theta^3 - 10 \theta^2 - 12\theta + 21}{26} \, .
\]
We perform a partial $2$-cover descent on $C_{6,-2}$ over $L$. The outcome
of this is that if $(X,Y)$ is a rational point then
there are non-zero $a \in \Q$ and $\alpha$, $\beta \in L$ such that
\[
X-\phi_1 = \mu_1 \cdot a \cdot \alpha^2,
\qquad X-\phi_2=\mu_2 \cdot a \cdot \beta^2.
\]
We note in passing that the only points on the
 curve $C_{6,-2}$ appear to be $\{ (-1,\pm 4)\}$, and 
that
\[
-1-\phi_1= 3 \mu_1 
\frac{(3 \theta^4 + 2 \theta^3 - 15 \theta^2 - 10\theta + 13)^2}{6^2},
\quad
-1-\phi_2= 3\mu_2 
 \frac{(\theta^3 + 2\theta^2 - 4\theta - 9)^2}{3^2} \, ,
\]
which provides a useful check on the correctness of our partial descent
implementation.

The ring of integers of $L$ is $\Z[\theta]$. The ideal
$2 \cdot \Z[\theta]$ factors as $\fp \fq^2$, where $\fp$ and $\fq$
are prime ideals. 
If $r \in \Q$ is non-zero,
then $\ord_\fq(r)=2 \ord_2(r)$; in particular $\ord_\fq(X)$
and $\ord_\fq(a)$ are even.
We note that $\ord_\fq(\mu_1)=1$. Thus $\ord_\fq(X-\phi_1)$
is odd. However, $\ord_\fq(\phi_1)=0$. This forces $\ord_2(X)=0$.
As $X=v/w$ and $v$, $w$ are coprime, we see that $v$ and $w$
are both odd. Now \eqref{eqn:case6} forces $2s^2 \equiv 2 t^2 \equiv
0 \pmod{4}$, contradicting the coprimality of $s$ and $t$. 

\begin{prop} \label{prop:f6}
  The only solutions in integers $(s,t,u)$ with $s$ and~$t$ coprime and satisfying
  $s \not\equiv t \bmod 2$ and $t \not\equiv 0 \bmod 3$ to the equation
  \[ u^5 = 3 s^4 + 6 s^2 t^2 - t^4 \]
  are $(s,t,u) = (0, \pm 1, -1)$.
This yields the trivial solution $(x,y,z)=(0,-1,-1)$ to
\eqref{eqn:main}. 
\end{prop}

\bigskip

This completes the proof of Theorem~\ref{thm:main}.


\end{document}